\documentclass{article}

\newcommand{\R}{I\!\! R}

\newcommand{\Z}{I\!\! Z}
\newcommand{\N}{I\!\! N}

\begin{document}

Tsemo Aristide

The International Center for Theoretical Physics

Strada Costiera, 11

Trieste, Italy.

email tsemo@ictp.trieste.it

\medskip

New Address:

 Tsemo Aristide,

3738, Avenue de Laval Appt 106

Montreal Canada H2X 3C9

email: tsemoaristide@hotmail.com

\bigskip

\centerline{\bf Closed similarity lorentzian affine manifolds.}

\bigskip

\centerline {\bf Abstract.}

\medskip

{\it A $Sim(n-1,1)$ affine manifold is a $n-$dimensional affine
manifold whose linear holonomy lies in the similarity Lorentzian
group but not in the Lorentzian group. In this paper, we show that
a compact $Sim(n-1,1)$ affine manifold is incomplete. Let $<,>_L$
be the Lorentz form, and $q$ the map on ${\R}^n$ defined by
$q(x)=<x,x>_L$. We show that for a compact radiant $Sim(n-1,1)$
affine manifold $M$ whose developing map is injective, if a connected component $C$ of
${\R}^n-q^{-1}(0)$ intersects the image of the universal cover of
$M$ by the developing map, then either $C$ or a connected
component of $C-H$, where $H$ is an hyperplane is contained in
this image.}

\bigskip

\centerline{\bf Introduction.}

\medskip

An $n-$dimensional affine manifold $M$, is an $n-$dimensional
differentiable manifold endowed with an atlas whose coordinate
changes are locally affine maps. The affine structure of $M$ pulls
back to its universal cover $\hat M$, and defines on it an affine
structure determined by a local diffeomorphism $D:\hat M
\rightarrow {\R}^n$, called the developing map. The developing map
gives rise to a representation $h:\pi_1(M)\rightarrow
Aff({\R}^n)$, called the holonomy of the affine manifold. Its
linear part $L(h)$, is called the linear holonomy of the affine
manifold. We will say that the affine manifold is complete, if and
only if the developing map is a diffeomorphism. An $n-$affine
manifold is said to be radiant if its holonomy fixes an element of
${\R}^n$.

We denote by $O(p,q)$, the subgroup of linear automorphisms of
${\R}^n$ which preserve a bilinear symmetric form of type $p,q$,
and by $Sim(p,q)$ the group generated by $O(p,q)$ and the
homotheties. An $O(p,q)$ affine manifold $M$ is an affine manifold
$M$ such that the image of its linear holonomy $L(h)$ is a
subgroup of
 $O(p,q)$.
 An $Sim(p,q)$ affine manifold $M$ is an affine
 manifold $M$ such that the image of its linear holonomy $L(h)$
is a subgroup of $Sim(p,q)$,  and contains an element which is not
in $O(p,q)$.

Let consider the flat riemmannian torus $T^n$, Bieberbach has
shown that closed  $O(n,0)$ affine manifolds are finitely covered
by $T^n$. Using the notion of discompacity, Yves Carri\`ere has
shown that closed $O(n-1,1)$ affine manifolds are complete. It is
obvious that a $Sim(n,0)$ affine manifold is incomplete, since an
element of its holonomy which doesn't lie in $O(n,0)$ fixes an
element of ${\R}^n$. There exist examples of complete $Sim(n-1,1)$
affine manifolds. Let's give one:

Endow ${\R}^n$  with its basis $(e_1...,e_n)$ and with the
lorentzian product defined by
$$
<e_i,e_i>_L=1; 0<i<n; <e_i,e_j>_L= 0; i\neq j; <e_n,e_n>_L=-1.
$$

We restrict this product to ${\R}^2$.
The affine map whose linear part is
$$
\pmatrix{2&0\cr 0&1}$$ in the basis  $(e_1+e_2,e_1-e_2)$, and
whose translation part is $e_1-e_2$ generates  a group which acts
properly and freely on ${\R}^2$.

The goal of this paper is to study closed $Sim(n-1,1)$ affine
manifolds. First we show:

\medskip

{\bf Theorem 1.} {\it A compact $Sim(n-1,1)$ affine manifold is
incomplete.}

\medskip

 After, using
the notion of  discompacity, we show

\medskip

{\bf Theorem 2.} {\it Let $M$ be a closed radiant $Sim(n-1,1)$
affine manifold whose developing map is injective,   if a connected component $C$ of
${\R}^n-q^{-1}(0)$ intersects $D(\hat M)$, then either $C$ is
contained in $D(\hat M)$ or a connected component of $C-H$, where
$H$ is an hyperplane.}

\bigskip

Interesting structures of $Sim(n-1,1)$ affine manifolds can be
constructed using the work of Goldman on projective structures on
surfaces  see [Gol]. For instance a $Sim(2,1)$ structure which
linear holonomy  is Zariski dense in $Gl(3,{\R})$ is given in this
paper.

\medskip

{\bf 1. Closed $Sim(n-1,1)$ affine manifolds are incomplete.}

\medskip

The main goal of this part is to show that a closed $Sim(n-1,1)$
affine manifold cannot be complete.

Let suppose that there exists a complete closed $Sim(n-1,1)$
affine manifold $M$; $M$ is the quotient of ${\R}^n$ by a subgroup
of affine transformations $\Gamma$, whose linear part is contained
in $Sim(n-1,1)$.

\medskip

{\bf Lemma 1.1.} {\it Let $\gamma$ be an element of $\Gamma$ whose
linear part has a determinant $<1$. Then there exists a basis
$(e_1,...,e_n)$ of ${\R}^n$ such that the linear part of $\gamma$
in this basis has the following form:
$$
\pmatrix{1&0&0\cr 0&{1\over {\lambda}^2}&0\cr 0& 0& {1\over
{\lambda}}B"}
$$
where $\lambda$ is a real number strictly superior to $1$ in
absolute value, and $B"$ is a matrix which preserves the
restriction of an euclidean product to  the sub vector space
generated by $e_3,...,e_n$.}

\medskip

{\bf Proof.}

We have supposed that the  determinant of the linear part
 $L(\gamma)$ of $\gamma$, is strictly inferior to $1$
 in absolute value. This implies that
 there exists a real number $\lambda>1$ such that
 $ \lambda L(\gamma)= L(\gamma)'$, where $L(\gamma)'$ is an element of
 $O(n-1,1)$. The linear map $L(\gamma)$ has $1$ as eigenvalue, since
 $\gamma$ acts freely. We deduce that $\lambda$ is an eigenvalue of
 $L(\gamma)'$. We remark that $L(\gamma)'$ has another eigenvalue
 $\alpha$
 which module is
 different from $1$ and the module of $\lambda$ since the absolute value
 of its determinant is $1$.
 If $\alpha$ is not a real number, then $\alpha$ and its complex conjugated
 $\bar\alpha$ are eigenvalues associated to the complex eigenvectors
 $u_1$ and $u_2$. In this case the restriction of $L(\gamma)'$ to
 the plane generated by $u_1+u_2$ and $i(u_1-u_2)$ is an euclidean similitude
 whose ratio is different from $1$. This is impossible since
 $L(\gamma)'$ lies in $O(n-1,1)$.
 Let $v_1$ and $v_2$ be the eigenvectors associated to $\lambda$ and
 $\alpha$ and $<,>_L$ the lorentzian product preserved by the linear
 holonomy. We have:
 $$
 <v_1,v_1>_L=<v_2,v_2>_L=0.$$

 We deduce that the restriction of $<,>_L$ to the plane $P$ generated by
  $v_1$ and $v_2$ is nondegenerate and has signature $(1,1)$. This implies that the restriction
  of $<,>_L$ to the orthogonal $W$ of $P$ with respect to
  itself is a scalar product. The restriction
  $B"$ of $L(\gamma)'$ to $W$ is
 an orthogonal linear map. We can suppose that its determinant is $1$.
 We deduce that
  $\alpha={1\over\lambda}$.

  \medskip

  Up to a change of origin, we can suppose that
  $\gamma(0)=(a_1,0...,0)$ where $a_1$ is a real number. The restriction
  $B$ of $L(\gamma)$ to the linear subspace generated by
  $(e_2,...,e_n)$ is strictly contracting.
  It is easy to show that the group generated by $\gamma$ is not
  cocompact, so $\Gamma$ contains another element $\gamma_1$ different
  from $\gamma$.

  \medskip

  {\bf Lemma 1.2.}
  {\it Let $C$ be the linear part of $\gamma_1$, then $C(e_1)=e_1+b$
  where $b$ lies in the linear subspace generated by
  $e_2,...,e_n$. }

  \medskip

  {\bf Proof.}

  Let $k$ be an element of ${\N}$.
  Consider the element $\gamma^k\gamma_1$. Its linear part has $1$
  as eigenvalue. The matrix
  of this linear part in the basis $(e_1,...,e_n)$ is $A^kC$, where
  $A$ is the matrix of the linear part of $\gamma$.
  Let $u_k$ be an eigenvector of $A^kC$ associated to $1$. We assume that
  the norm of $u_k$ with respect to the euclidean scalar product
  defined $<e_i,e_j>=\delta_{ij}$  is $1$.
  Let $u_k=(u^1_k,u^2_k)$. We have  $C(u_k)=(v^1_k,v^2_k)$,
   where $u^1_k$ and $v^1_k$ are elements of ${\R}$, and
  $u^2_k$ and $v^2_k$ are elements
   of the vector space  $F$ generated by $e_2,...,e_n$.
  We have $v^1_k=u^1_k$ since $A^k(e_1)=e_1$, and $A^k$ preserves $F$.
  Since $B$ is strictly contracting the norm of $A^k(0,v^2_k)$ goes to $0$
  with respect to the euclidean norm. So $u_k$ goes to $e_1$ and
  $C(e_1)$, which is the limit of $C(u_k)=A^{-k}(u_k)$ is $e_1+b$ where $b$ is
  an element of $F$.

  \medskip

  {\bf First proof of the theorem 1.}

  Let $c$ be the translational part of $\gamma_1$ in the basis
  $(e_1,...,e_n)$.
  Put $c=(c_1,...,c_n)$.
  We have $\gamma^k\circ\gamma_1\circ\gamma^{-k}(0)=
  (c_1,B^k(-ka_1b+(c_2,...,c_n)))$. Since $B$ is contracting and the action
  of $\Gamma$ is proper, we deduce that $(c_2,...,c_n)=b=0$.
    This implies that
    $\gamma^m\gamma_1^n(0)=(na_1+mc_1)e_1$. Since the action of $\Gamma$
    on ${\R}^n$ is proper and free we deduce that the subgroup 
    $\{n,m\in{\Z}, na_1+mc_1\}$ is discrete. We deduce from this fact that
    there exist $p, q$ in ${\Z}$ such that $pa_1+qc_1=0$. This implies that
    $\gamma^p={\gamma_1}^{-q}$.

  Let $K$ be a fundamental domain of the action of $\Gamma$, recall that it is
  a compact such that for each $x\in{\R}^n$, there exists an element $\gamma_x$
  of $\Gamma$ such that $\gamma_x(x)\in K$, and for each element $\gamma$ of
  $\Gamma$, $\gamma(K^0)\cap K^0$ is empty, where $K^0$ is the interior of $K$.
  Let $\mid\mid$, $\mid\mid$ be a norm associated to  $<,>$, the  scalar product for which $(e_1,...,e_n)$ is an orthonormal
  basis. There exists a real $A>0$ and an integer $l$ such that
   $\mid\mid K\mid\mid< A$ and $l\mid a_1\mid>A$. Consider the element 
   $u=(0,y_0)$ of ${\R}^n$ where $y_0$ is an element in the vector subspace
   generated by $(e_3,..,e_n)$ if $n\geq 3$, otherwise $y_0$ is an element 
   of $Vect(e_2)$ such that $\mid\mid B^l(y_0)\mid\mid>A$.
   There must exist an element $\gamma_0$ in $\Gamma$ such that
    $\gamma_0(u)$ is in $K$. We have shown that there exist elements
    $p, q$ in ${\Z}$ such that ${\gamma_0}^p=\gamma^q$.
    Denote by $D$ the restriction of $\gamma_0$ to $Vect(e_2,..,e_n)$,
    and $d_1$ the real number such that $\gamma_0(0)=d_1e_1$.
    We have $D^p=B^q$ and $pd_1=qa_1$. Since $\mid\mid y_0\mid\mid>A$,
     $D$ is contractant and $q>p$. It results from lemma 1.1 that the restriction
     of $D$ and $B$ to $Vect(e_3,..,e_n)$ are similarities of respective
     ratio $r_d$ and $r_b$ (If the dimension is $2$, we consider the 
     restriction of $D$ and $B$ to $Vect(e_2)$ which is a similarity).
     We have $(r_d)^p=(r_b)^q$. Moreover $\mid\mid B^l(y_0)\mid\mid=
     (r_b)^l\mid\mid y_0\mid\mid>A$, and $A>\mid\mid D(y_0)\mid\mid=
     r_d\mid\mid y_0\mid\mid$. We deduce that $(r_b)^l>r_d$. This implies
     that $(r_b)^{lq}>(r_d)^{q}$, which is equivalent to saying that
     $(r_b)^{lq}=(r_d)^{lp}>(r_d)^q$. This implies that $q>lp$. 
     
    $p$ and $q$ are elements of ${\N}$. We have
    that $\gamma_0(u)=(d_1,D(y_0))\in K$. This implies that
    $A>\mid d_1\mid$. But we have $qA>lpA>lp\mid d_1\mid=lq\mid a_1\mid$.
    Which is a contradiction since we have supposed that $l\mid a_1\mid>A$.

  \medskip
  
  We can also deduce the proof of Theorem 1, from this deep result 
  
  \medskip
  
  {\bf Theorem [F-G-H].}
  
  {\it Let $M$ a compact affine manifold, which is the quotient of 
  ${\R}^n$ by the group $\Gamma$ which acts properly and freely; then
  $\Gamma$ does not preserve proper affine subspaces.}
  
  \medskip
  
  {\bf Proof.}
  Up to a finite cover, we can assume that $M$ is oriented
   We denote by $(C(M),d)$ the
  simplicial complex we define the simplicial homology of $M$. We can lift the simplicial
  decomposition to ${\R}^n$, and thus lift the simplicial complex $(C(M),d)$
  to the simplicial complex $(C({\R}^n,d')$. This last complex has a structure
  of a ${\Z}\Gamma$ module. Its homology is trivial since ${\R}^n$ is 
  contractible. We deduce that it is a resolution of $\Gamma$. We can use this
  resolution to calculate the real cohomology of $\Gamma$. The cohomology
  obtained is also the real cohomology of $M$. The De Rham theorem implies that 
  ${H^n}_{DR}(M,{\R})=H^n(M,{\R})=H^n(\Gamma,{\R})={\R}$ since we have supposed that $M$ is 
  oriented. Suppose that $\Gamma$ preserves an $l-$affine subspace $F$ of ${\R}^n$.
  We denote by $N$ the quotient of $F$ by $\Gamma$. We have also that
  $H_{DR}^n(N,{\R})=H^n(\Gamma,{\R})={\R}$. Since $N$ is compact, its implies
  that $l\geq n$, we deduce that $l=n$.

  \medskip
  
  {\bf Second proof of theorem 1.}
  
  Let $c$ be the translational part of $\gamma_1$ as in the first proof,
  we can remark that in the basis $(e_1,...,e_n)$, we have
  $c=(c_1,0,..0)$ this implies that $\Gamma$ preserves the line ${\R}e_1$.
  This fact contradicts the previous theorem.
  
  \medskip

{\bf Remark.}

 The both proofs of theorem 1 are related. While proving theorem 1, we
 have shown that for every element $\gamma_1$ in $\Gamma$, there exists
 $p$, and $q$ in ${\Z}$ such that $\gamma_1^p=\gamma_q$ this implies that
 the quotient of $\Gamma$ by the group generated by $\Gamma$ is a torsion
 group. We deduce  that the real cohomological dimension of $\Gamma$ is
 one which is contrary to the fact that up to a finite cover 
 $H^n(M,{\R})=H^n(\Gamma,{\R})={\R}$.

 In contrast to the $Sim(n,0)$ affine manifolds (See [Fr] theorem 1),
  there exist compact
 $Sim(n-1,1)$ affine manifolds which are not radiant. Here is an example.

Endow ${\R}^2$ with the Lorentzian product $(,)$ such that
$(e_1,e_1)=(e_2,e_2)=0$, and $(e_1,e_2)=1$.

 Consider the subgroup $\Gamma$ of $Aff({\R}^2)$ generated by the following
 transformations:

 $$\gamma_1(x,y)=(x+1,y),$$
 $$\gamma_2(x+y)=(x,2y),$$
 the quotient of ${\R}\times ({\R}-\{0\})$ by $\Gamma$ is a compact
 $Sim(n-1,1)$ affine manifold.

 \bigskip

  {\bf 2. On the universal cover of compact $Sim(n-1,1)$ affine manifolds.}

  \medskip

  In this part we are going to find properties of the universal cover  of
  a closed radiant $Sim(n-1,1)$ affine manifold.  We use the notion
  of discompacity defined by Carri\`ere [Car] 2.2.1. Let us recall it.

  We consider in ${\R}^n$ the unit ball $B_n$. The euclidean metric
  induces on closed subsets of ${\R}^n$ the Hausdorff distance.
  Let $G$ be a subgroup of $Gl(n,{\R})$, and $(g_p)_{p\in{\N}}$ a
  sequence of elements of $G$.
  The limit of the family $(g_p(B_n)\cap B_n)_{p\in {\N}}$  converges in
  $B_n$. It
  is a degenerated ellipsoid  (see [Car]). The codimension
   of this ellipsoid is  the discompacity $d$, of the
    family $(g_p)_{p\in{\N}}$, the discompacity of the group with respect to
   the euclidean metric is the smallest $d$.

   Obviously we cannot use the notion of discompacity in this form since
   the linear holonomy of our manifold
   may contain homotheties. Denote $q:{\R}^n\rightarrow {\R}$
   $x\rightarrow <x,x>_L$. We can define in ${\R}^n-q^{-1}(0)$ the metric
   $$
   (u,v)\longrightarrow <,u,v>'={<u,v>_{euc}\over { q(x)}}
   $$
   where $u$ and $v$ are vectors of the tangent space at $x$ and
   $<,>_{euc}$ is the euclidean scalar product.

   \medskip

  {\bf  Theorem 2.1.}

  {\it  Let $\hat x$ be an element of $\hat M$, $u$ and $v$, elements of
  $T_{\hat x}{\hat M}$, such that the geodesics
  $c_1: [0,1]\rightarrow \hat M$, $t \rightarrow exp_{\hat x}(tu)$, and the one
  $c_2:[0,1]\rightarrow \hat M$, $t \rightarrow  exp_{\hat x}(tv)$
  are defined. Suppose that
   the elements $exp_{\hat x}(u)$
  and $exp_{\hat x}(v)$ can't be joined by a geodesic, but for every $t,t'<1$,
  there is a geodesic between $exp_{\hat x}(tu)$ and $exp_{\hat x}(t'v)$.
  Let $c:[0,1]\rightarrow {\R}^n$, $t\rightarrow exp_{D(exp(u))}(tw)$ be
  the geodesic between $D(exp_{\hat x}(u))$ and $D(exp_{\hat
  x}(v))$, and
  let $U_{\hat x}$ be the domain of definition of $exp_{\hat x}$.
 Consider the element $t_0\in [0,1]$ such that for every $t<t_0$,
 $exp_{D(exp_{\hat x}(u))}(tw)$ is an element of
  $D(exp_{\hat x}(U_{\hat x}))$
 but not is $exp_{D(exp_{\hat x}(u))}(t_0w) =y$. Then $y$ is an element of
 $q^{-1}(0)$.}

\medskip

{\bf Proof.}

There is a geodesic $\hat c_3:[0,1[\rightarrow \hat M$,
$t\rightarrow exp_{\hat x}(tb)$ such that $y$ is an element of the
adherence of $D(\hat c_3([0,1[))$ and such that $D(\hat c_3)$ is
contained in the convex hull of $D(\hat c_1)$ and $D(\hat c_2)$,
where $\hat c_1$, and $\hat c_2$ are geodesics of $\hat M$
respectively above $c_1$ and $c_2$. Set $p(\hat x)=x$, the image
$c_3$ of $p(\hat c_3)$ is a maximal incomplete geodesic of $M$.
Since $M$ is compact, there exists an element $z$ of $M$ such that
the geodesic $c_3$ is recurrent in an affine chart $U$ which
contains $z$. We deduce as Carri\`ere,  the existence of a family
of ellipsoids $s_p$ of ${\R}^n$ whose centers are elements of
$D(\hat c_3)$, such that for each $p,p'$, there is an element
$\gamma_{p,p'}$ of the holonomy such that
$\gamma_{p,p'}(s_p)=s_{p'}$ and the centers $x_p$ of $s_p$ goes to
$y$.

Suppose that $y$ is not an element of $q^{-1}(0)$.

Let $z_p$ be an element of an ellipsoid $s_p$, and $u_p$, $v_p$
two vectors in its tangent space. Put
$\gamma_{p,p'}=\lambda_{p,p'}g_{p,p'}$ where $g_{p,p'}$ is an
element of $O(n-1,1)$. We have:
$$
{<\gamma_{p,p'}(u_p),\gamma_{p,p'}(v_p)>_{euc}\over
q(\gamma_{p,p'}(x))}= {<g_{p,p'}(u_p),g_{p,p'}(v_p)>_{euc} \over
q(x)}.
$$
Since the holonomy of $M$ is supposed to be radiant.

 The metrics $<,>_{euc}$ and $<,>'$
 are equivalent in a neighborhood of $y$
 since $q(y)$ is different from $0$.
 We know that the discompacity of the family of $g_p$
 in respect to the riemmannian metric $<,>_{euc}$ is $1$.
 The family of ellipsoids $s_p$ goes to an ellipsoid, or a
 codimension $1$ degenerated ellipsoid centered in $y$.
 We conclude as in Carri\`ere that $y$  must be an element of
$D(exp_{\hat x}(U_{\hat x}))$. This is not possible, so $q(y)=0.$

\medskip

A similar result is given in [Gol].

\medskip

{\bf Corollary 2.2.} {\it Let $M$ be a compact radiant
$Sim(n-1,1)$ affine manifold, let $\hat x$, $u$, and $v$  be
respectively elements of $\hat M$ and  $T_{\hat x}\hat M$, such
that $exp_{\hat x}(u)$ and $exp_{\hat x}(v)$ are defined.
 If the convex hull  $E$ of $(D(\hat x), D(exp_{\hat x}(u)),
 D(exp_{\hat x}(v))$ is contained in a connected component
 of ${\R}^n-q^{-1}(0)$, then it is contained in
 $D(exp_{\hat x}(U_{\hat x}))$.}

\medskip

{\bf Proof.}

Suppose that $E$ is not contained in $D(exp_{\hat x}(U_{\hat
x}))$. Let $y$ and $z$ be two elements of $E\cap D(exp_{\hat
x}(U_{\hat x}))$ such that $y=D(exp_{\hat x}(u_1))$,
$z=D(exp_{\hat x}(u_2))$, and for every $t_1, t_2<1$, $exp_{\hat
x}$ is defined on the convex hull of $0, tu_1, tu_2$, but
  the elements $exp_{\hat x}(u_1)$ and $exp_{\hat x}(u_2)$
   cannot be joined by
a geodesic. Consider the geodesic $c:[0,1]\rightarrow {\R}^n$,
$t\rightarrow exp_y(tw)$ between $y$ and $z$. There  exists a real
number $0<t_0<1$, such that for $0<t<t_0$, $exp_y(tw)$ lies in
 $D(exp_{\hat x}(U_{\hat x}))$, but
  not $exp_y(t_0w)$. We deduce from the theorem 2.1. that $exp_y(t_0w)$
  must lie in $q^{-1}(0)$. This is contrary to the hypothesis.

\medskip

More generally, we can determine, how the boundary of the image of
the developing map of a compact radiant $Sim(n-1,1)$ affine
manifold is: more precisely, we have the following proposition
which implies theorem 2:

\medskip

{\bf Proposition 2.3.} {\it Let $M$ be a compact radiant
$Sim(n-1,1)$ affine manifold whose developing map is injective,  the boundary
 of $D(\hat M)$, is contained in  the union of $q^{-1}(0)$
 and an hyperplane.}

 \medskip

 {\bf Proof.}

 As in [Car] p. 625, one can remark that elements of the boundary of
 $D(\hat M)$ which are not elements of $q^{-1}(0)$
 are limits of $(\gamma_ne)_{n\in {\N}}$, where $\gamma_n$ is an element of
 the holonomy and $e$ is an ellipsoid. We conclude that those
 elements are contained in at most two hyperplane $H_1$, $H_2$. The case of
 two hyperplane is impossible, since those hyperplane are stable
 by the holonomy, the affine function $\alpha$ such that
 $\alpha(H_1)=0$ and $\alpha(H_2)=1$, will be invariant by the holonomy
 and so define a differentiable function on $M$ without maximal.
 (It is the same argument used in [Car]).

 \medskip

{\bf Proposition 2.4.} {\it Let $M$ be a compact radiant affine
manifold,  if the image of the developing map is a convex set
contained in an open set of ${\R}^n-q^{-1}(0)$, then the
developing map is injective.}

\medskip

{\bf Proof.}

Let $\hat x$ be an element of $\hat M$. For every elements, $u$
and $v$ of $U_{\hat x}$, the convex hull of $D(\hat x)$,
$y=D(exp_{\hat x}(u))$ and $z=D(exp_{\hat x}(v))$ is a subset of
$D(\hat M)\cap ({\R}^n-q^{-1}(0)).$ We deduce from the corollary
$2.2$ that $y$ and $z$ are elements of $D(U_{\hat x})$. This
implies that $U_{\hat x}$ is a convex set. We can conclude by
using [Kos].

\medskip

A particular case of the situation of corollary 2.4 is the
following: endow  a compact oriented surface $S$ of genus $>2$,
with an hyperbolic structure, and consider $q$ the Lorentzian form
defined on ${\R}^3$ by $q(x_1,x_2,x_3)={x_1}^2+{x_2}^2-{x_3}^2$.
The hyperbolic structure can be defined by a representation of the
fundamental group of $S$, $\pi_1(S)\rightarrow O(2,1)$ such that
the quotient of $H=q^{-1}(-1)$ by $\pi_1(S)$ is $S$. The quotient
of $W=\{x: q(x)<0, x_3>0\}$ by the group generated by $\pi_1(S)$,
and a homothetie of ratio $0<\lambda<1$, is a compact $Sim(n-1,1)$
affine manifold whose universal cover is $W$.

More generally  we have

\medskip

{\bf Corollary 2.5.} {\it Let $M$ be a radiant compact affine
manifold such that the image of its developing  map is contained
in $W=\{x: q(x)<0, z>0\}$, $M$ is the quotient of a connected
component of  $W-H$ by a discrete group of $Sim(n-1,1)$, where $H$
is an hyperplane of ${\R}^n$.}

\medskip

{\bf Proof.}
We remark that the interior of
a connected component of $W-H$ is convex.
This implies the image of the developing map is a convex set. The result
follows using 2.3 and 2.4.

\medskip

Let $M$ be a compact radiant $Sim(n-1,1)$ affine manifold, the
foliation $D(\hat {\cal F}_q)$ of ${\R}^n-\{0\}$ whose leaves are
 the sub manifolds defined by $q=$ constant is invariant by the holonomy
 of $M$. Its pull back on $\hat M$ defines a foliation
 $\hat{\cal F}_q$ of $\hat M$, which gives rise to a foliation
 ${\cal F}_q$ of $M$. If $D(\hat N)=D(\hat M)\cap q^{-1}(0)$
  is not empty, then $N=p(D^{-1}(D(\hat N)))$ is a compact submanifold
  of $M$. Note that the $1-$parameter group $\phi_t$ generated by
  the radiant vector field $X_0$, preserves the foliation
  ${\cal F}_q$, and is transverse to all the leaves but not to the
  connected components of $N$.

  \medskip

  {\bf Proposition 2.6.}
  {\it If $N$ is empty, then $M$ is the total space of a bundle over
  $S^1$.}

  \medskip

  {\bf Proof.}

  If $N$ is empty, then  $\phi_t$ is transverse
  to the foliation ${\cal F}_q$. This implies that this foliation is
  a Lie foliation. We conclude using [God] Corollary 2.6 p. 154.

  \medskip

  {\bf Remark.}

  Recall that a riemannian foliation on a manifold $M$, is a foliation
  ${\cal F}$ of $M$ such that there exists a riemanian metric of $M$ which
  projects locally along the leaves of $M$, this is equivalent to saying that
  locally the distance between the leaves is well-defined, or the foliation
  is defined by locally submersions $M\rightarrow N$ which transitions functions
  preserve a riemannian metric of $N$.
  
  Suppose that the image of the developing map is included in
the upper cone  $C=\{ x/ q(x)<0\,x_n>0\}$, one can define a map
$f:C\rightarrow q^{-1}(-1)$ such that $f(x)$ is the element of $C$ colinear to
$x$ such that $q(f(x))=-1$. Let $(,)_L$ be the flat Lorentzian product of
${\R}^n$, The restriction of the lorentzian product
defined by $(u,v)_x={{(u,v)_L}_x\over (x,x)_L}$ to
$q^{-1}(-1)$ is an hyperbolic metric. This endows the radial flow 
$\hat \phi_t$ of $C$ with a transverse riemannian structure. Since the holonomy
of the manifold is included in ${\R}O(n-1,1)$, the riemannian structure
of $\hat \phi_t$ pushes forward to a riemannian structure of its radial
flow $\phi_t$.
In fact the flow $\phi_t$ is also a transversally $(O(n-1,1),q^{-1}(-1))$
homogeneous foliation where $O(n-1,1)$. This eanables
one to define the global holonomy $h_{\phi}$ of $\phi_t$ which is a 
representation $h_{\phi}:\pi_1(M)\rightarrow O(n-1,1)$. It assigns to any
element $\gamma$ of $\pi_1(M)$ (identified in this case as a subgroup
of ${\R}O(n-1,1))$, the element $t\gamma (t\in {\R})$ such that $t\gamma$ is an
element of $O(n-1,1)$ and it preserves $C$.

The riemannian foliation  have been intensively studied, It has been shown by
Molino that the adherence of a leaf of a riemannan foliation defined on a
compact manifold is a submanifold. Carriere and Carron have shown that in the
 case of riemannian flows, the adherence of leaves are torus.
 
 More precisely, for a riemanian foliation, let denote by $(M_1,{\cal F}_1)$
 the bundle of transverse orthogonal frames of $M$ endowed with the 
 pulls back ${\cal F}_1$ of ${\cal F}$, one can define the sheaf of local
 vector fields of $M_1$ which commute with the global foliated vector fields
 of ${\cal F}_1$. This sheaf pushes forward to a sheaf $C(M,{\cal F})$ of
 $M$. It is called the commuting sheaf of ${\cal F}$. The second structure
 theorem of riemannian foliations says that the adherences of the leaves of
 ${\cal F}$ are orbits of the pseudogroup defined by $C(M,{\cal F})$.
 In the case of a riemannian flow, the Lie algebra $C(M,{\cal F})$ is 
 commutative.
 
 In the other hand  consider the adherence $L$ of the image of $h_{\phi}$ in $O(n-1,1)$.
 If the image of $h_{\phi}$ is discrete, it implies that the orbit of $\phi_t$
 are closed, the holonomy of a leaf of $\phi_t$ is finite, this implies that
 up to a finite cover we can consider that the foliation $\phi_t$ does not
 have holonomy. This finite cover of $M$ is the total space of a bundle 
 whose typical fiber is an hyperbolic manifold. This manifold is the quotient
 of $q^{-1}(-1)$ by the image of $h_{\phi}$.
 
 If $L$ is not discrete, then its lie algebra $l$ is isomorphic to 
 $C(M,{\cal F})$ see [W], we deduce that this Lie algebra is commutative. This 
 imply that the connected component $L_0$ of $L$ is a non trivial commutative
 group. Since $\pi_1(M)$ normalizes $L_0$, using [G-K] 1.3 one can conclude that
 $\pi_1(M)$ is a solvable group. We have shown:
 
 \medskip
 
 {\bf Theorem.}
 
 {\it Suppose that the image of the developing map of a $Sim(n-1,1)$
 compact radiant affine manifold is contained in $C=\{x\in {\R}^n, 0>q(x), and x_n>0\}$,
 then if $\pi_1(M)$ is not solvable, the leaves of the radiant flow are 
 compact.}
 
 \medskip
 
 {\bf Remark.}
 
 The previous result is a particular case of a result due to Epstein for
 transversely hyperbolic foliation. Moreover if Epstein shows that if the
 leaves of the radiant flow are not compact, then the dimension of $M$ is $3$ or
 $4$, and $M$ is a the quotient of a solvable group $G$ endowed with a 
 left symmetric structure.

\bigskip

\centerline{\bf References.}

\medskip

[Car] Carri\`ere, Y. Autour de la conjecture de L. Markus sur les
vari\'et\'es affines. Invent. math. 95, (1989) 615-628.

[Car-Car] Caron, P. Carriere, Y. Flots transversalement de Lie ${\R}^n$,
flots de Lie transversalement minimaux

[E] Epstein, D. Transversally hyperbolic 1-dimensional foliation Asterisque
116.

[Fr] Fried, D. Closed similarity manifolds.
Comment. Math. Helvetici 55, (1980) 576-582.

[F-G-H] Fried, D. Goldman, W. Hirsch, M. Affine manifolds with
nilpotent holonomy. Comment. Math. Helvetici 56, (1981) 487-523.

[God] Godbillon, C. Feuilletages, \'etudes g\'eom\'etriques, I.
Progress. in Math. vol 98. (1991).

[Gol] Goldman, W.  Projective structures with Fuchsian holonomy.
J. Differential Geometry 25, (1987) 297-326.

[G-K] Goldman, W. Kamishima, Y. The fundamental group of a compact
flat Lorentz space form is virtually polycyclic. J. Diff. Geom. 19
(1984) 233-240.

[Kos] Koszul, J-L. Vari\'et\'es localement plates et convexit\'e.
Osaka J. Math. 2, (1965) 285-290.

[L] Leslie, J. A remark on the group of automorphisms of a foliation
having dense leaf. Jour. of Diff. Geom., 7 (1972) 597-601

[Mo] Molino, P. Riemannian foliation, Progress in Mathematics (1988)

[W] Wolak, R. Foliated $G-$structures and Riemannian foliations, 
Manuscripta. Math. 66 (1989), 45-59.

\end{document}